\documentclass[12pt]{article}
\usepackage{latexsym,amsfonts,amssymb}
\usepackage[dvips]{color}
\usepackage{colordvi,multicol}
\usepackage{amsmath}

\def \cal{\mathcal}

\newtheorem{thm}{Theorem}[section]

\newtheorem{defn}[thm]{Definition}

\newtheorem{exa}[thm]{Example}

\date{}
\begin{document}
\title{\bf Periodic solutions of hybrid jump diffusion processes}
\author{}

\maketitle

\centerline{Xiao-Xia Guo}
\centerline{\small
School of Applied Mathematics}
\centerline{\small Shanxi University of Finance and Economics}
\centerline{\small  Taiyuan, 030006, China}
\centerline{\small  and}
\centerline{\small
School of Mathematics and Information
Sciences}
\centerline{\small Guangzhou University}
\centerline{\small  Guangzhou, 510006, China}
\centerline{\small E-mail: xxguo91@163.com}

\vskip 1cm \centerline{Wei Sun}
\centerline{\small  Department of
Mathematics and Statistics}
\centerline{\small  Concordia University}
\centerline{\small Montreal, H3G 1M8, Canada}
\centerline{\small
E-mail: wei.sun@concordia.ca}

\begin{abstract}

\noindent In this paper, we investigate periodic solutions of
regime-switching jump diffusions. We first show the well-posedness
of solutions to SDEs corresponding to the hybrid system. Then, we
derive the strong Feller property and the irreducibility of the
associated time-inhomogeneous semigroups. Finally, we establish
the existence and uniqueness of periodic solutions. Concrete
examples are presented to illustrate the results.
\end{abstract}

\noindent  {\it MSC:} 60J25; 60H10; 60J75; 34C25

\noindent  {\it Keywords:} Hybrid system; regime-switching jump diffusion;
periodic solution; strong Feller property; irreducibility.


\section{Introduction}
A hybrid system is a dynamical system whose evolution depends on
both continuous and discrete variables. The study of hybrid
systems is becoming more and more  important in different research
areas such as biology, ecosystems, wireless communications, signal
processing, engineering and mathematical finance. In the recent
years, lots of progress has been made on a class of hybrid systems
called regime-switching jump diffusion processes. This model
consists of two component processes $(X(t),\Lambda(t))$ with
$X(t)$ and $\Lambda(t)$ being of continuous and discrete states,
respectively. If $\Lambda(t)$ is a Markov chain that is
independent of $X(t)$, then we have the model of Markov-switching
jump diffusions; whereas if $\Lambda(t)$ depends on $X(t)$, then
we have the model of regime-switching jump diffusions. We refer
the reader to the  monographs \cite{fam,YZ} for comprehensive
studies of hybrid switching diffusions and their applications.

The model of regime-switching jump diffusions provides more
flexibility in applications but also requires careful examination
of the dependence  between the continuous and discrete components.
In the past decade, various properties of these hybrid systems
have been thoroughly studied and many remarkable results have been
obtained. For example, Chen et al. obtained maximum principles and
Harnack inequalities in \cite{ACFZ0} and discussed the recurrence
and ergodicity in \cite{ACFZ}. Xi \cite{ACFA}, Xi and Yin
\cite{BBC} investigated asymptotic properties of the model. Nguyen
and Yin \cite{NY1,NY2}, Shao \cite{AAE}, Xi and Zhu \cite{FFRTA},
Xi, Yin and Zhu \cite{zvaqw} considered the model whose switching
component has a countably infinite state space. Note that most of
the existing literatures have focused on the time-homogeneous
case. In this paper, we will study time-inhomogeneous
regime-switching jump diffusions and investigate their periodic
solutions.

Periodic solutions play an  important role in the study of
stochastic dynamical systems. Here we just list some previous
works that are closely related to this paper. Khasminskii
\cite{a8} systematically developed the theory of periodic
solutions to random systems modelled by stochastic differential
equations (SDEs). In \cite{k1}, Zhang et al.  investigated
periodic solutions of SDEs driven by L\'{e}vy processes. Hu and Xu
\cite{a12} obtained the existence and uniqueness of periodic
solutions for stochastic functional differential equations by
proving the global attractivity of solutions.  In our recent work
\cite{fazxc}, we established the ergodicity and uniqueness of
periodic solutions for jump diffusions by proving  the strong
Feller property and the irreducibility of the associated
time-inhomogeneous semigroups. Note that none of the above
mentioned works discusses periodic solutions of hybrid systems.
Our present paper seems to be the first one discussing periodic
solutions of regime-switching jump diffusions whose switching
components can have countably infinite state spaces.

The rest of this paper  is organized as follows. In Section 2, we
will show the well-posedness of solutions to SDEs corresponding to
the hybrid system  (\ref{aa}) and (\ref{a3}) (see below). The
unique strong solution is obtained by representing the switching
component as a stochastic integral with respect to a Poisson
random measure (see  (\ref{dts}) below) and by using an
interlacing procedure. Different from \cite[Proposition 2.1]{BBC},
we only assume that the coefficient functions satisfy the local
Lipschitz condition. Also, we remove a key assumption of
\cite{AAE,zvaqw,FFRTA}, which requires that the transition matrix
of the switching component is H\"{o}lder continuous. In Section 3,
we will establish the existence and uniqueness of periodic
solutions. To this end, we derive the strong Feller property and
the irreducibility of the associated time-inhomogeneous
semigroups. Our assumption $\mathbf{Q_1}$ (see below) is weaker
than \cite[(A1)]{AAE}, \cite[Assumption 5.1]{zvaqw} and
\cite[Assumption 4.3]{FFRTA}. Finally, we will give two examples
in Section 4 to demonstrate our main results.

\section{Existence and uniqueness of regime-switching jump diffusion processes}

Let $(\Omega,\mathcal{F},\{\mathcal{F}_t\}_{t\geq0},\mathbb{P})$
be a complete probability space with  filtration
$\{\mathcal{F}_t\}_{t\geq0}$ satisfying the usual conditions
(i.e., it is increasing, right continuous and $\mathcal{F}_0$
contains all $\mathbb{P}$-null sets). Suppose that $k,l, m\in
\mathbb{N}$ with $k\geq m$. Denote by $\mathbb{R}_+$ the set of
all non-negative real numbers. Let $\{B(t)\}_{t\ge 0}$ be a
$k$-dimensional standard Brownian motion and $N$ be an independent
Poisson random measure (corresponding to a stationary point
process $p(t)$) on $\mathbb{R}_+\times(\mathbb{R}^l-\{0\})$. The
compensator $\widetilde{N}$ of $N$ is given by
$\widetilde{N}(\mathrm{d}t,\mathrm{d}u)=N(\mathrm{d}t,\mathrm{d}u)-\nu(\mathrm{d}u)\mathrm{d}t$,
where $\nu(\cdot)$ is a L\'{e}vy measure satisfying
$\int_{\mathbb{R}^l-\{0\}}(1\wedge|u|^2)\nu(\mathrm{d}u)<\infty$.
Let $\mathcal{S}=\{1,2,\dots\}$.

Throughout this paper, we fix a $\theta>0$. Let $(X(t), \Lambda(t))$ be a pair of right continuous processes
with left-hand limits on $\mathbb{R}^m\times \mathcal{S}$.
The first component $X(t)$ satisfies the
following SDE:
\begin{eqnarray}\label{aa}
dX(t)&=&b(t,X(t),\Lambda(t))\mathrm{d}t+\sigma(t,X(t),\Lambda(t))\mathrm{d}B(t)
+{\int_{\{|u|<1\}}}H(t,X(t-),\Lambda(t-),u)\widetilde{N}
(\mathrm{d}t,\mathrm{d}u)\nonumber\\
&& +{\int_{\{|u|\geq 1\}}}G(t,X(t-),\Lambda(t-),u)N
(\mathrm{d}t,\mathrm{d}u).
\end{eqnarray}
We assume that the coefficient functions
$b(t,x,i):[0,\infty)\times\mathbb{R}^m\times\mathcal{S}\rightarrow\mathbb{R}^m $,
$\sigma(t,x,i):[0,\infty)\times\mathbb{R}^m\times\mathcal{S}\rightarrow
\mathbb{R}^{m\times k}$,
$H(t,x,i,u):[0,\infty)\times\mathbb{R}^m\times\mathcal{S}\times\mathbb{R}^l
\rightarrow\mathbb{R}^m$ and
$G(t,x,i,u):[0,\infty)\times\mathbb{R}^m\times\mathcal{S}\times\mathbb{R}^l
\rightarrow\mathbb{R}^m$ are all Borel measurable and satisfy
\begin{eqnarray*}
&&b(t+\theta,x,i)=b(t,x,i),\ \ \ \ \sigma(t+\theta,x,i)=\sigma(t,x,i),\nonumber\\
&&H(t+\theta,x,i,u)=H(t,x,i,u),\ \ \ \ G(t+\theta,x,i,u)=G(t,x,i,u)
\end{eqnarray*}
for any $t\ge 0$, $x\in\mathbb{R}^m$, $i\in\mathcal{S}$ and $u\in\mathbb{R}^l-\{0\}$.
The second component $\Lambda(t)$ has the  state space $\mathcal{S}$
such that
\begin{eqnarray}\label{a3}
\mathbb{P}\{\Lambda(t+\Delta)=j|\Lambda(t)=i,X(t)=x\}
=\left\{\begin{array}{ll}
          q_{ij}(x)\Delta +o(\Delta),& i\neq j, \\
          1+ q_{ij}(x)\Delta +o(\Delta),& i=j,
        \end{array}\right.
\end{eqnarray}
as $\Delta\downarrow0$. Hereafter $q_{ij}(x)$ is a Borel
measurable function on $\mathbb{R}^m $ for $i,j\in{\cal S}$ such
that $q_{ij}(x)\geq0$  for all $i,j\in\mathcal{S}$ with $i\neq j$
and $\sum_{j\in\mathcal{S}}q_{ij}(x)=0$ for all $i\in\mathcal{S}$
and $x\in\mathbb{R}^m$. We make the following assumption on the
matrix $Q=(q_{ij}(x))$.

\noindent $\mathbf{Q_1}$)
$$\mathop{\sum}\limits_{j=1}^\infty\mathop{\sup}\limits_{i\in\mathcal{S}\setminus\{j\},\,x\in\mathbb{R}^m}q_{ij}(x) <\infty.$$

We point out that $\Lambda(t)$ can be represented as a stochastic integral with respect to a Poisson random measure.
To this end, we construct a family of intervals $\{\Delta_{ij}(x):i,j\in\mathcal{S}, i\neq j\}$ on the positive real line as follows:
for $x\in\mathbb{R}^m$ and $i,j\in\mathcal{S}$  with $j\neq i$, set
\begin{eqnarray*}
  \Delta_{i1}(x) &=& [0,q_{i1}(x)), \\
   \Delta_{i2}(x)&=& [q_{i1}(x),q_{i1}(x)+q_{i2}(x)), \nonumber\\
   \vdots&&  \nonumber\\
  \Delta_{ij}(x) &=& \left[\sum_{s=1}^{j-1}q_{is}(x),\sum_{s=1}^{j}q_{is}(x)\right),\nonumber\\
  \vdots&&  \nonumber\\
\end{eqnarray*}
If $q_{ij}(x)=0$ for $j\neq i$, we set $\Delta_{ij}(x)=\emptyset$.
Define
$$a(j):=\sup_{i\in\mathcal{S}\setminus\{j\},\,x\in\mathbb{R}^m}q_{ij}(x),\
\ \ \  L=\sum_{j=1}^{\infty}a(j).
 $$By the assumption $\mathbf{Q_1}$, we get $L<\infty$. Then, each value of the interval $\Delta_{ij}(x)$ is bounded by $L$.

 We define a function
$h : \mathbb{R}^m\times\mathcal{S}\times [0,L] \rightarrow\mathbb{R}$ by
$h(x,i,r)=\sum_{j\in\mathcal{S}}(j-i)1_{\Delta_{ij}(x)}(r)$, that is,
$$ h(x,i,r) = \left\{\begin{array}{ll}
      j-i, &{\rm if}\ r\in\Delta_{ij}(x), \\
                            0,& \mathrm{otherwise}. \\
                          \end{array}\right.
$$
Then, $\Lambda(t)$ can be modeled by
\begin{eqnarray}\label{dts}
  d\Lambda(t) &=& \int_{[0,L]}h(X(t-),\Lambda(t-),r)N_1(dt,dr),
\end{eqnarray}
where $N_1(dt, dr)$ is a Poisson random measure (corresponding to
a stationary point process $p_1(t)$ which is adapted to
$\mathcal{F}_t$) with characteristic measure ${\cal L}(dr)$, the
Lebesgue measure on $[0,L]$. The Poisson random measure
$N_1(\cdot,\cdot)$ is assumed to be independent of the Brownian
motion $B(\cdot)$ and the Poisson random measure $N(\cdot,
\cdot)$. Therefore, the processes $(X(t), \Lambda(t))$ modeled by
(\ref{aa}) and (\ref{a3}) can be thought of as a solution to the
hybrid system (\ref{aa}) and (\ref{dts}).

We use $|x|$ to denote the Euclidean norm of a vector $x$, use
$A^T$ to denote the transpose of a matrix $A$, and use
$|A|:=\sqrt{{\rm trace}(A^TA)}$ to denote the trace norm of $A$.
Let
$C^{1,2}(\mathbb{R}_+\times\mathbb{R}^m\times\mathcal{S};\mathbb{R})$
be the space of all real-valued functions $f(t,x,i)$ on
$\mathbb{R}_+\times\mathbb{R}^m\times\mathcal{S}$ which are
continuously differentiable with respect to $t$ and twice
continuously differentiable with respect to $x$. For $f\in
C^{1,2}(\mathbb{R}_+\times\mathbb{R}^m\times\mathcal{S};\mathbb{R})$,
define
\begin{eqnarray*}
{\cal A}f(t,x,i):={\cal L_i}f(t,x,i)+Q(x)f(t,x,\cdot)(i),
\end{eqnarray*}
with
\begin{eqnarray*}
&&{\cal L_i}f(\cdot,\cdot,i)(t,x)\\&:=&f_t(t,x,i)+\langle f_x(t,x,i), b(t,x,i)\rangle+
\frac{1}{2}\mathrm{ trace}(\sigma^T(t,x,i)f_{xx}(t,x,i)\sigma(t,x,i))\nonumber\\
&&+\int_{\{|u|<1\}}[f(t,x+H(t,x,i,u),i)-f(t,x,i)-
\langle f_x(t,x,i), H(t,x,i,u)\rangle]\nu(\mathrm{d}u)\nonumber\\
&&+\int_{\{|u|\ge 1\}}[f(t,x+G(t,x,u,i),i)-f(t,x,i)]\nu(\mathrm{d}u),
\end{eqnarray*}
and
\begin{eqnarray*}
Q(x)f(t,x,\cdot)(i):=\sum_{j\in\mathcal{S}}q_{ij}(x)[f(t,x,j)-f(t,x,i)],
\end{eqnarray*}
where $f_t=\frac{\partial f}{\partial t}$, $f_x=\nabla_x
f=(\frac{\partial f}{\partial x_1},\dots,\frac{\partial
f}{\partial x_m})$ and $f_{xx}=(\frac{\partial^2 f}{\partial
x_i\partial x_j})_{m\times m}$. Define a metric $d(\cdot,\cdot)$
on $\mathbb{R}^m\times\mathcal{S}$ by
$d((x,i),(y,j))=|x-y|+|i-j|$. A Borel measurable function $g$ on
$[0,\infty)$ is said to be locally integrable, denoted by $g\in
L_{ loc}^1([0,\infty);\mathbb{R})$, if
$$
\int_0^{\tau}|g(x)|dx<\infty,\ \ \forall \tau>0.
$$

We make the following assumptions.

\noindent$\mathbf{A_1}$)\ \ For each $i\in{\cal S}$,
\begin{eqnarray*}
b(\cdot,0,i)~,\sigma(\cdot,0,i)\in
{L}^2([0,\theta);\mathbb{R}^m),
~~\int_{\{|u|<1\}}|H(\cdot,0,i,u)|^2\nu(\mathrm{d}u)\in
{L}^1([0,\theta);\mathbb{R}^m).
\end{eqnarray*}
For each $n\in\mathbb{N}$, there exists
$L_{n}(t)\in L^1([0,\theta);\mathbb{R}_+)$ such that for any $t\in
[0,\theta)$, $i\in\mathcal{S}$ and $x,y\in\mathbb{R}^m$ with $|x|\vee|y|\leq n$,
\begin{eqnarray*}
 &&   | b(t,x,i)-b(t,y,i)|^2\leq
  L_{n}(t)|x-y|^2,\ \ \ \
|\sigma(t,x,i)-\sigma(t,y,i)|^2\leq
  L_{n}(t)|x-y|^2, \nonumber \\
 &&   \int_{\{|u|<1\}}|H(t,x,i,u)-H(t,y,i,u)|^2\nu(\mathrm{d}u)\leq
L_{n}(t)|x-y|^2.
\end{eqnarray*}

\noindent $\mathbf{B_1}$)\ \ There exist $V_1\in
C^{1,2}(\mathbb{R}_+\times\mathbb{R}^m;\mathbb{R}_+)$ and
$g_{1}\in L_{ loc}^1([0,\infty);\mathbb{R})$ such that
\begin{eqnarray}\label{da}
 \lim_{|x|\rightarrow\infty}
\left[\inf_{t\in[0,\infty)}V_1(t,x)\right]=\infty,
\end{eqnarray}
and for any $t\ge 0$, $i\in\mathcal{S}$ and $x\in\mathbb{R}^m$
\begin{eqnarray}\label{dc}
{\cal L}_iV_1(t,x)\leq g_{1}(t).
\end{eqnarray}

\begin{thm}\label{lem-1}
Suppose that assumptions $\mathbf{A_1}$, $\mathbf{B_1}$ and $\mathbf{Q_1}$ hold.
Then, the hybrid system given by (\ref{aa}) and (\ref{dts}) has a unique strong solution $(X(t),\Lambda(t))$ with initial value $(X(0),\Lambda(0))=(x,i)$.
Moreover, $\mathbb{P}\{T_\infty=\infty\}=1$, where $T_\infty=\lim_{n\rightarrow\infty}T_n$ and $T_n=\inf\{t\geq 0:|X(t)|\vee\Lambda(t)\geq n\}$.
\end{thm}

\noindent {\bf Proof.} Let
$(x,i)\in\mathbb{R}^m\times\mathcal{S}$. By \cite[Theorem
2.2]{fazxc}, under assumptions $\mathbf{A_1}$ and $\mathbf{B_1}$,
there exists a unique non-explosive strong solution $X^{(i)}(t)$
to the following SDE:
\begin{eqnarray}\label{dd}
dX^{(i)}(t)&=&b(t,X^{(i)}(t),i)\mathrm{d}t+\sigma(t,X^{(i)}(t),i)\mathrm{d}B(t)
+{\int_{\{|u|<1\}}}H(t,X^{(i)}(t-),i,u)\widetilde{N}
(\mathrm{d}t,\mathrm{d}u)\nonumber\\
&& +{\int_{\{|u|\geq 1\}}}G(t,X^{(i)}(t-),i,u)N
(\mathrm{d}t,\mathrm{d}u)
\end{eqnarray}
with initial value $X^{(i)}(0)=x$. Let $\sigma_1 < \sigma_2 <
\cdots < \sigma_n < \cdots$ be the set of all jump points
 of the
stationary point process $p_1(t)$ corresponding to the Poisson
random measure $N_1(dt, dr)$. Since ${\cal L}(\cdot)$ is a finite
measure on $[0, L]$, $\lim_{n\rightarrow\infty}\sigma_n=+\infty$
almost surely.

 We now construct the processes $(X(t),\Lambda(t))$. Define
\begin{equation}\label{er1}
(X(t),\Lambda(t))=(X^{(i)}(t),i),\ \ t \in[0,
 \sigma_1).
\end{equation}
Let
$$
\Lambda(\sigma_1)=i+
\sum_{j\in\mathcal{S}}(j-i)1_{\Delta_{ij}(X^{(i)}(\sigma_1-))}(p_1(\sigma_1)).
$$
Define
\begin{equation}\label{er2}
(X(t),\Lambda(t))=(X^{(i)}(\sigma_1),\Lambda(\sigma_1)),\ \
t=\sigma_1.
\end{equation}
Note that (\ref{aa}) is
equivalent to (\ref{dd}) on the time interval $[0, \sigma_1)$. Hence the above defined processes
$(X(t),\Lambda(t))$ provide the unique strong solution to the
hybrid system (\ref{aa}) and (\ref{dts}) on $[0, \sigma_1)$.

Let
$$
\widetilde{B}(t)=B(t+\sigma_1)-B(t), \ \
\widetilde{p}(t)=p(t+\sigma_1),\ \
\widetilde{p_1}(t)=p_1(t+\sigma_1).$$ Set
\begin{eqnarray*}
  \begin{array}{l}
     (\widetilde{X}(t),\widetilde{\Lambda}(t)) =(X^{(\Lambda(\sigma_1))}(t),\Lambda(\sigma_1)),\ \ t\in[0,\sigma_2-\sigma_1), \\
     \widetilde{X}(\sigma_2-\sigma_1)=X^{(\Lambda(\sigma_1))}(\sigma_2-\sigma_1),\\
     \widetilde{\Lambda}(\sigma_2-\sigma_1)=\Lambda(\sigma_1)+\sum\limits_{j\in\mathcal{S}}
  (j-\Lambda(\sigma_1))\mathrm{1}_{\widetilde{A}(j)}(\widetilde{p}_1(\sigma_2-\sigma_1)),
  \end{array}
\end{eqnarray*}
where
$$\widetilde{A}(j)=\Delta_{\Lambda(\sigma_1),j}(X^{(\Lambda(\sigma_1))}((\sigma_2-\sigma_1)-)).$$
Then, we define
$$(X(t),\Lambda(t))=(\widetilde{X}(t-\sigma_1),\widetilde{\Lambda}(t-\sigma_1)), \ \ t\in[\sigma_1,\sigma_2],$$
which together with  (\ref{er1})  and (\ref{er2}) gives the unique
strong solution on the time interval $[0, \sigma_2]$. Continuing
this procedure inductively, we define $(X(t), \Lambda(t))$ on the
time interval $[0, \sigma_n]$ for each $n$. Therefore, $(X(t),
\Lambda(t))$ is the unique strong solution to the hybrid system
(\ref{aa}) and (\ref{dts}).

Next, we show that $(X(t),\Lambda(t))$ is non-explosive. Define
$$\tau_0=0,\ \ \tau_p = \inf\{t
>\tau_{p-1} : \Lambda(t)\neq \Lambda(\tau_{p-1})\},\ \ \tau
= \lim_{p\rightarrow\infty}\tau_p. $$
 Then, $\tau$ is the first
instant prior to which $\Lambda(t)$ has infinitely many switches.
Note that $\{\tau_p\}_{p\geq1}$ is a subsequence of
$\{\sigma_n\}_{n\geq1}$ since $\Lambda(t)$ can have jumps only at
the time sequence $\{\sigma_{n}\}_{n\geq1}$. Define $q_i(x):=-q_{ii}(x)$ for $x\in \mathbb{R}^m$ and $i\in {\cal S}$. Based on the interlacing construction of $(X(t), \Lambda(t))$, we get
\begin{eqnarray}\label{pop}\mathbb{P}\{\tau_{p+1}-\tau_{p}>t\}=
\mathbb{E}\left[\mathrm{exp}\left\{-\int_0^tq_{\Lambda(\tau_{p})}(X^{\Lambda(\tau_{p})}(s))ds\right\}\right]
,\end{eqnarray}
and
\begin{eqnarray*}
\mathbb{P}\{\Lambda(\tau_{p+1})=j|\mathcal{F}_{\tau_{p+1}-}\}&=&
\frac{q_{\Lambda(\tau_{p}),j}(X(\tau_{p+1}-))}{q_{\Lambda(\tau_{p})}(X(\tau_{p+1}-))}(1-\delta_{\Lambda(\tau_{p}),j})1_{\{q_{\Lambda(\tau_{p})}(X(\tau_{p+1}-))>0\}}\\
&&+\delta_{\Lambda(\tau_{p}),j}1_{\{q_{\Lambda(\tau_{p})}(X(\tau_{p+1}-))=0\}}.
\end{eqnarray*}

By the assumption $\mathbf{Q_1}$, for any $ i\in\mathcal{S}$, we have
\begin{eqnarray*}
&&q_i(x)=-q_{ii}(x)=\sum_{j\in\mathcal{S}\setminus\{i\}}q_{ij}(x)\leq \sum_{j=1}^\infty a(j)=L.
\end{eqnarray*}
Then, (\ref{pop}) implies that $\mathbb{P}\{\tau_{p+1}-\tau_{p}>t\}\geq e^{-Lt}$ for all $p\in\mathbb{N}$
and $t > 0$. Hence
\begin{eqnarray*}
\mathbb{P}\{\tau_\infty=\infty\} &\geq& \mathbb{P}\{\{\tau_{p+1}-\tau_{p}>t\}\ \ \mathrm{i.o.}\}\nonumber\\
&=&  \mathbb{P}\left\{\bigcap_{r=1}^{\infty}\bigcup_{s=r}^{\infty}\{\tau_{s+1}-\tau_{s}>t\}\right\}\nonumber\\
&=& \lim_{r\rightarrow\infty}\mathbb{P}\left\{\bigcup_{s=r}^{\infty}\{\tau_{s+1}-\tau_{s}>t\}\right\} \nonumber\\
&\geq& \limsup_{r\rightarrow\infty}\mathbb{P}\{\tau_{r+1}-\tau_{r}>t\}\\
&\geq& e^{-Lt}.
\end{eqnarray*}
Letting $t\downarrow0$, we get $\mathbb{P}\{\tau_\infty=\infty \}=1$.
Therefore,
$(X(t),\Lambda(t))$ is the unique strong solution to the hybrid system (\ref{aa}) and (\ref{dts}) with initial value
$(x,i)$ for all $t\in[0,\infty)$.

Finally, we show that $\mathbb{P}\{T_\infty=\infty\}=1$.
If this is not true, then
there exist $\varepsilon>0$ and $T'\in (0,\infty)$ such that
$$\mathbb{P}\{T_\infty\leq T'\}>2\varepsilon.$$
Hence we can find a sufficiently large integer $n_0$ such that
\begin{eqnarray}\label{bz}
 \mathbb{P}\{T_n\leq T'\}>\varepsilon,\ \ \forall n\geq n_0.
\end{eqnarray}
Define an increasing sequence by
$$\rho_1=1,\ \ \rho_{n}=\min\bigg\{r\in \mathbb{N}:r\geq \rho_{n-1}+2,\ \
\sum_{j=r}^\infty a(j)\le\frac{1}{2^{n}}\bigg\}.$$
Set $f(j)=n$ for $j\in [\rho_n,\rho_{n+1})$ and $n\in\mathbb{N}$. Then, $f$ is a non-decreasing function on ${\cal S}$ satisfying $f(j)\rightarrow\infty$ as $j\rightarrow\infty$ and
$$\sum_{j=1}^{\infty}a(j)f(j)
=\sum_{n=1}^{\infty}\sum_{j\in[\rho_n,\rho_{n+1})}a(j)f(j)\leq\sum_{n=1}^{\infty}\frac{n}{2^n}<\infty.$$

Define $V(t,y,k):=V_1(t,y)+f(k)$, where $V_1$ is given
in the assumption $\mathbf{B_1}$.
By It\^{o}'s formula, (\ref{dc}) and the assumption $\mathbf{Q_1}$, we get
\begin{eqnarray*}
&&\mathbb{E}_{0,x,i}[V(t\wedge T_n,X(t\wedge T_n),\Lambda(t\wedge T_n))]\\
 &=&V(0,X(0),\Lambda(0))+\mathbb{E}_{0,x,i}\left[\int_0^{t\wedge T_n}{\cal A}V(s,X(s),\Lambda(s))ds\right] \\
 &=& V(0,x,i)+\mathbb{E}_{0,x,i}\left[\int_0^{t\wedge T_n}{\cal L}_{\Lambda(s)}V_1(s,X(s))ds\right]\\
  &&+\mathbb{E}_{0,x,i}\left[\int_0^{t\wedge T_n}\sum_{j\in\mathcal{S}}q_{{\Lambda(s),j}}(X(s))[f(j)-f(\Lambda(s))]ds\right]\\
&\leq& V(0,x,i)+\mathbb{E}_{0,x,i}\left[\int_0^{t\wedge T_n}g_1(s)ds\right] +\mathbb{E}_{0,x,i}\left[\int_0^{t\wedge T_n}\sum_{j=1}^\infty a(j)f(j)ds\right].
\end{eqnarray*}
Then, we have
\begin{eqnarray*}
  \mathbb{E}_{0,x,i}[V(T'\wedge T_n,X(T'\wedge T_n),\Lambda(T'\wedge T_n))]
  \leq V(0,x,i)+\int_0^{T'}g_1(s)ds +T'\sum_{j=1}^\infty a(j)f(j).
\end{eqnarray*}
Set
$$
M:=V(0,x,i)+\int_0^{T'}g_1(s)ds +T'\sum_{j=1}^\infty a(j)f(j)<\infty.
$$
Then,
\begin{eqnarray}\label{bx}
\mathbb{E}_{0,x,i}[1_{\{T_n\leq T'\}}V(T_n,X(T_n),\Lambda(T_n))]\leq M.
\end{eqnarray}

Define
$$\mu(n)=\inf\{V(t,y,j):(t,y,j)\in\mathbb{R}_+\times\mathbb{R}^m\times\mathcal{S},~|y|\vee j\geq n\}.$$
Then, $\lim_{n\rightarrow\infty}\mu(n)=\infty$ by
(\ref{da}) and our choice  of $f$. However, by (\ref{bz}) and (\ref{bx}), it follows that
$$\varepsilon\mu(n)<\mu(n)\mathbb{P}\{T_n\leq T'\}\leq M.$$
We have arrived at a contradiction.
Therefore,
\begin{equation}\label{beta}
T_\infty=\infty\ \ \ \ a.s..
\end{equation}
The proof is complete.
\hfill\fbox

\section{Periodic solutions of regime-switching jump diffusion processes}\label{sec2}\setcounter{equation}{0}

In this section, we will study the existence and uniqueness of periodic
solutions of the hybrid system given by (\ref{aa}) and (\ref{dts}). Denote by
$\mathcal{B}(\mathbb{R}^m\times\mathcal{S})$ the Borel $\sigma$-algebra of
$\mathbb{R}^m\times\mathcal{S}$ and denote by $B_b(\mathbb{R}^m\times\mathcal{S})$ (respectively
$C_b(\mathbb{R}^m\times\mathcal{S})$)  the space of all real-valued bounded Borel
functions (respectively continuous and bounded functions) on
$\mathbb{R}^m\times\mathcal{S}$.
Suppose that $\{Y(t),t\geq0\}$ is a Markov process on $\mathbb{R}^m\times\mathcal{S}$. We define its transition probability function by
$$
P(s,y,t,A)=\mathbb{P}\{Y(t)\in
A|Y(s)=y\},\ \ y\in\mathbb{R}^m\times\mathcal{S},\,0\leq s<t,
$$
and define the corresponding Markovian semigroup of linear operators $\{P_{s,t}\}$ on $B_b(\mathbb{R}^m\times\mathcal{S})$ by
 $$(P_{s,t}f)(y):=\mathbb{E}_{s,y}[f(Y(t))]:
 =\int_{\mathbb{R}^m\times\mathcal{S}}f(z)P(s,y,t,dz),\ \ y\in\mathbb{R}^m\times\mathcal{S},\,0\leq s<t.$$

\begin{defn}
(i) A Markov process $\{Y(t),t\geq0\}$ on $\mathbb{R}^m\times\mathcal{S}$ is said to be
$\theta$-periodic if for any $n\in\mathbb{N}$ and any $0\le
t_1<t_2<\dots<t_n$, the joint distribution of the random variables
$Y(t_1+k\theta),Y(t_2+k\theta),\ldots, Y(t_n+k\theta)$ is
independent of $k$ for $k\in \mathbb{N}\cup\{0\}$.
A Markovian transition semigroup $\{P_{s,t}\}$  is said to be
$\theta$-periodic if $P(s,y,t,A)=P(s+\theta,y,t+\theta,A)$
for any $0\le s< t$, $y\in\mathbb{R}^m\times\mathcal{S}$ and $A\in
\mathcal{B}(\mathbb{R}^m\times\mathcal{S})$.

(ii) A family of probability measures $\{\mu_s, s\ge 0\}$ is said to be $\theta$-periodic with respect to
the Markov semigroup $\{P_{s,t}\}$ of $Y(t)$ if
$$
\mu_s(A)=\int_{\mathbb{R}^m\times\mathcal{S}}
P(s,y,s+\theta,A)\mu_s(dy),\ \ \forall s\ge 0,\, A \in
\mathcal{B}(\mathbb{R}^m\times\mathcal{S}).
$$

(iii) A stochastic process $\{Y(t),t\geq0\}$ is said to be a
$\theta$-periodic solution of the system given by (\ref{aa})  and
(\ref{dts}) if it is a solution of the system given by (\ref{aa})
and (\ref{dts}) and is $\theta$-periodic.
\end{defn}

\begin{defn}
Let $0\le s_0<t_0<\infty$. A Markovian transition semigroup
$\{P_{s,t}\}$ is said to be regular at $(s_0,t_0)$ if all
transition probability measures $P(s_0,y,t_0,\cdot)$, $y\in\mathbb{R}^m\times\mathcal{S}$,
are mutually equivalent. $\{P_{s,t}\}$ is said to be strongly Feller at
$(s_0,t_0)$ if $P_{s_0,t_0}f\in
C_b(\mathbb{R}^m\times\mathcal{S})$ for any $f\in B_b(\mathbb{R}^m\times\mathcal{S})$. $\{P_{s,t}\}$ is said to be irreducible  at
$(s_0,t_0)$ if $P(s_0,y,t_0,A)>0$ for any $y\in \mathbb{R}^m\times\mathcal{S}$ and
any non empty open subset $A$ of $\mathbb{R}^m\times\mathcal{S}$. $\{P_{s,t}\}$ is
said to be regular, strongly Feller, irreducible if it is
regular, strongly Feller, irreducible at any $(s_0,t_0)$,
respectively.
\end{defn}

\subsection{Strong Feller and irreducible properties of time-inhomogeneous semigroups}

By Theorem \ref{lem-1}, under assumptions $\mathbf{A_1}$, $\mathbf{B_1}$ and $\mathbf{Q_1}$, there exists a unique strong solution $(X(t),\Lambda(t))$ to
the hybrid system (\ref{aa}) and (\ref{dts}). By the interlacing structure, one finds that $(X(t),\Lambda(t))$ is a Markov process on $\mathbb{R}^m\times\mathcal{S}$.
In this subsection, we will show that the transition semigroup $\{P_{s,t}\}$ of $(X(t),\Lambda(t))$ is strongly Feller and irreducible.

We make the following assumptions.

\noindent $\mathbf{A_2}$)\ \
For each $i\in\mathcal{S}$,
\begin{eqnarray*}
&&b(\cdot,0,i)\in
{L}^2([0,\theta);\mathbb{R}^m),\, \sigma(\cdot,0,i)\in
{L}^{\infty}([0,\theta);\mathbb{R}^m), \nonumber\\
&&\int_{\{|u|<1\}}|H(\cdot,0,i,u)|^2\nu(\mathrm{d}u)\in
{L}^1([0,\theta);\mathbb{R}^m).
\end{eqnarray*}
For each $n\in\mathbb{N}$, there exists
$L_n\in L^{\infty}([0,\theta);\mathbb{R}_+)$ such that for any $t\in
[0,\theta)$, $i\in\mathcal{S}$ and $x,y\in\mathbb{R}^m$ with $|x|\vee|y|\leq n$,
\begin{eqnarray}\label{qww}
 &&   | b(t,x,i)-b(t,y,i)|^2\leq
  L_n(t)|x-y|^2,\ \ \ \
|\sigma(t,x,i)-\sigma(t,y,i)|^2\leq
  L_n(t)|x-y|^2,  \nonumber\\
 &&   \int_{\{|u|<1\}}|H(t,x,i,u)-H(t,y,i,u)|^2\nu(\mathrm{d}u)\leq
L_n(t)|x-y|^2.
\end{eqnarray}

\noindent $\mathbf{A_3}$)\ \ For each $i\in \mathcal{S}$, $t\in[0,\theta)$ and $x\in\mathbb{R}^m$,
$Q(t,x,i):=\sigma(t,x,i)\sigma^T(t,x,i)$ is invertible and \begin{equation}\label{Q}
\sup_{|x|\leq
n,\,t\in[0,\theta)}| Q^{-1}(t,x,i)|<\infty,\ \ \forall n\in \mathbb{N},\ i\in \mathcal{S}.
\end{equation}

\noindent $\mathbf{B_2}$)\ \  There exists $V_2\in
C^{1,2}(\mathbb{R}_+\times\mathbb{R}^m;\mathbb{R}_+)$ such that
\begin{eqnarray*}
 \lim_{|x|\rightarrow\infty}
\left[\inf_{t\in[0,\infty)}V_2(t,x)\right]=\infty,
\end{eqnarray*}
and
$$
\sup_{x\in\mathbb{R}^m,\,t\in[0,\infty)}{\cal L}_iV_2(t,x)<\infty,\ \ \forall i\in\mathcal{S}.
$$

Obviously, $\mathbf{A_2}$ implies $\mathbf{A_1}$ and $\mathbf{B_2}$
implies $\mathbf{B_1}$.

\begin{thm}\label{thm-4} Suppose that assumptions $\mathbf{A_2}$, $\mathbf{A_3}$, $\mathbf{B_2}$ and $\mathbf{Q_1}$ hold.
Then the transition semigroup $\{P_{s,t}\}$ of $(X(t),\Lambda(t))$ is strongly Feller.
\end{thm}
{\bf Proof.} Denote the transition probability function of $(X(t),\Lambda(t))$ by $\{P(s,(x,i),t,B\times\{j\}):0\le s<t,(x,i)\in\mathbb{R}^m\times\mathcal{S},
B\in\mathcal{B}(\mathbb{R}^m),j\in{\cal S}\}$. For $(x,i)\in\mathbb{R}^m\times\mathbb{S}$, let $X^{(i)}(t)$ be defined by (\ref{dd}).
We kill the process $X^{(i)}(t)$ with rate $q_i(\cdot)$ and obtain a subprocess $\widetilde{X}^{(i)}(t)$ with generator $\mathcal{L}+q_{ii}$. Then,
\begin{eqnarray*}
  \mathbb{E}[f(\widetilde{X}^{(i)}(t))] &=& \mathbb{E}\left[f(X^{(i)}(t))\mathrm{exp}\left\{\int_0^tq_{ii}({X}^{(i)}(u))du\right\}\right],\ \ f\in\mathcal{B}_b(\mathbb{R}^m).
\end{eqnarray*}

Let $\widetilde{P}^{(i)}(s,x,\cdot)$
be the transition probability function of $\widetilde{X}^{(i)}(t)$.
Then, for $0\le s<t$, $B\in\mathcal{B}(\mathbb{R}^m)$ and $j\in{\cal S}$, we have
\begin{eqnarray*}
  P(s,(x,i),t,B\times\{j\}) &=&\delta_{ij}\widetilde{P}^{(i)}(s,x,t,B)\\
  &&+\int_s^t\int_{\mathbb{R}^m}P(t_1,(x_1,j_1),t,B\times\{j\})
  \sum_{j_1\in\mathcal{S}\setminus\{i\}}q_{ij_1}(x_1)\widetilde{P}^{(i)}(s,x,t_1,dx_1)dt_1.
\end{eqnarray*}
Repeating this procedure, we get
\begin{eqnarray*}
  P(s,(x,i),t,B\times\{j\})=\delta_{ij}\widetilde{P}^{(i)}(s,x,t,B) + \sum_{k=1}^{n}\Psi_k+U_n,
\end{eqnarray*}
where
$$
  \begin{array}{ll}
    \Psi_k= \displaystyle\mathop{\int\cdots\int}\limits_{s<t_1<\cdots<t_k<t}
   \sum_{j_0=i,j_1\in\mathcal{S}
   \setminus\{j_0\},j_2\in\mathcal{S}
   \setminus\{j_1\},\dots,j_k\in\mathcal{S}
   \setminus\{j_{k-1}\},j_k=j}\int_{\mathbb{R}^m}\cdots\int_{\mathbb{R}^m}
   \widetilde{P}^{(i)}(s,x,t_1,dx_1)\\
   ~~~~~~~~~\times q_{ij_1}(x_1)\widetilde{P}^{(j_1)}(t_1,x_1,t_2,dx_2)\cdots
   q_{j_{k-1}j_k}(x_{k})\widetilde{P}^{(j_k)}(t_{k},x_k,t,B)dt_1dt_2\cdots dt_k,
\end{array}$$
and
$$\begin{array}{ll}
     U_n =\displaystyle\mathop{\int\cdots\int}\limits_{s<t_1<\cdots<t_{n+1}<t}
   \sum_{j_0=i,j_1\in\mathcal{S}
   \setminus\{j_0\},j_2\in\mathcal{S}
   \setminus\{j_1\},\dots,j_{n+1}\in\mathcal{S}
   \setminus\{j_{n}\}}\int_{\mathbb{R}^m}\dots\int_{\mathbb{R}^m}
   \widetilde{P}^{(i)}(s,x,t_1,dx_1)\\
   ~~~~~~~~~\times q_{ij_1}(x_1)\widetilde{P}^{(j_1)}(t_1,x_1,t_2,dx_2)\cdots
   q_{j_{n}j_{n+1}}(x_{n+1})P(t_{n+1},(x_{n+1},j_{n+1}),t,B\times\{j\})dt_1dt_2\cdots dt_{n+1}.
  \end{array}
$$

By the assumption $\mathbf{Q_1}$, we find that $U_n$ does not
exceed $\frac{[(t-s)L]^{n+1}}{(n+1)!}$. Hence
\begin{eqnarray}\label{babao}
  P(s,(x,i),t,B\times\{j\})=\delta_{ij}\widetilde{P}^{(i)}(s,x,t,B) + \sum_{n=1}^{\infty}\Psi_n.
\end{eqnarray}
By \cite[Theorem 3.8]{fazxc} and assumptions $\mathbf{A_2}$,
$\mathbf{A_3}$, $\mathbf{B_2}$, we know that the transition
semigroup of $X^{(i)}(t)$ is strongly Feller. Following the
argument of \cite[Lemma 4.5]{zvaqw}, we can show that the
transition semigroup of $\widetilde{X}^{(i)}(t)$ is also strongly
Feller. Then, $\widetilde{P}^{(i)}(s,x,t,B)$ and $\Psi_n$, $n\in \mathbb{N}$, are all continuous with
respect to $x$. Note that
$\mathcal{S}$ is equipped with a discrete metric. Then, the
left-hand side of (\ref{babao}) is lower semi-continuous with
respect to $(x,i)$. Therefore, the transition semigroup
$\{P_{s,t}\}$ of $(X(t),\Lambda(t))$ is strongly Feller by
\cite[Proposition 6.1.1]{AAQ}.\hfill\fbox

Now we consider the irreducibility of the transition semigroup
$\{P_{s,t}\}$ of $(X(t),\Lambda(t))$. Denote by
$B_{b,loc}(\mathbb{R}_+)$ and
$B_{b,loc}([0,\infty)\times\mathbb{R}^m;\mathbb{R}^m)$ the sets of
all locally bounded Borel measurable functions on $\mathbb{R}_+$
and maps from $[0,\infty)\times\mathbb{R}^m$ to $\mathbb{R}^m$,
respectively. Let $f$ be a function on
$[0,\infty)\times\mathbb{R}^m$. For ${\rho}>0$, we define
$$f^{\bullet {\rho}}(t,x)=f(t,{\rho} x),\ \ t\ge
0,x\in\mathbb{R}^m.
$$
We make
the following assumptions.

\noindent $\mathbf{B_3})$\ \  There exists $V_3\in
C^{1,2}([0,\infty)\times\mathbb{R}^m;\mathbb{R}_+)$ satisfying the
following conditions:

  \noindent (i)
$$ \lim_{|x|\rightarrow\infty}
\left[\inf_{t\in[0,\infty)}V_3(t,x)\right]=\infty.
$$

 \noindent (ii) For any $\rho\ge 1$, there exist $q_{\rho}\in
 B_{b,loc}(\mathbb{R}_+)$ and $W_{\rho}(t,x)\in B_{b,loc}([0,\infty)\times\mathbb{R}^m;\mathbb{R}^m)$ satisfying for each $n\in\mathbb{N}$
there exists $R_n\in L_{loc}^1([0,\infty);\mathbb{R}_+)$ such that
for any $t\in [0,\infty)$ and $x,y\in\mathbb{R}^m$ with
$|x|\vee|y|\leq n$,
$$
   |W_{{\rho}}(t,x)-W_{{\rho}}(t,y)|^2\leq
  R_n(t)|x-y|^2,
  $$
and for $t\ge 0$, $x\in\mathbb{R}^m$ and $i\in{\cal S}$,
\begin{eqnarray*}
{\cal L}_iV^{\bullet {\rho}}_3(t,x)\leq q_{{\rho}}(t),
\end{eqnarray*}
\begin{eqnarray*}
V^{\bullet {\rho}}_3(t,x)\le\langle W_{{\rho}},\nabla_x V^{\bullet
{\rho}}_3(t,x)\rangle.
\end{eqnarray*}

\noindent $\mathbf{Q_2})$\ \ For any distinct $i,j
\in\mathcal{S}$, there exist $j_1,\dots, j_r \in \mathcal{S}$ with
$j_i \neq j_{i+1}$, $j_1=i$ and $j_r=j$ such that the set
$\{x\in\mathbb{R}^m : q_{j_{i}j_{i+1}}(x)>0\}$ has positive
Lebesgue measure for $i = 1,\ldots, r$.

\begin{thm}\label{thm-6} Suppose that assumptions $\mathbf{A_1}$, $\mathbf{A_3}$, $\mathbf{B_1}$, $\mathbf{B_3}$, $\mathbf{Q_1}$ and $\mathbf{Q_2}$ hold.
Then, the transition semigroup $\{P_{s,t}\}$ of $(X(t),\Lambda(t))$ is irreducible.
\end{thm}
{\bf Proof.} Denote the transition probability function of
$(X(t),\Lambda(t))$  by $\{P(s,(x,i),t,B\times\{j\}):0\le
s<t,(x,i)\in\mathbb{R}^m\times\mathcal{S},B\in\mathcal{B}(\mathbb{R}^m),j\in{\cal
S}\}$. As shown in the proof Theorem \ref{thm-4}, the assumption
$\mathbf{Q_1}$ implies that (\ref{babao}) holds. By \cite[Theorem
3.3.9]{fazxc}, we know that for any $i\in\mathcal{S}$, the
transition semigroup of $X^{(i)}(t)$ is irreducible under
assumptions $\mathbf{A_1}$, $\mathbf{A_3}$, $\mathbf{B_1}$ and
$\mathbf{B_3}$. Further, we find that the subprocess
$\widetilde{X}^{(i)}(t)$ is irreducible by the assumption
$\mathbf{Q_1}$. This together with the assumption $\mathbf{Q_2}$
implies that the right-hand side of (\ref{babao}) is positive
whenever $B$ is a non empty open set of $\mathbb{R}^m$. Then,
$P(s,(x,i),t,B\times\{j\})>0$. Since $B$ is arbitrary,  the
transition semigroup $\{P_{s,t}\}$ of $(X(t),\Lambda(t))$ is
irreducible.\hfill\fbox

\subsection{Existence and uniqueness of periodic solutions}
To prove the existence and uniqueness of periodic solutions to the hybrid system (\ref{aa}) and (\ref{dts}), we make the following assumptions.

\noindent $\mathbf{Q_3})$
$$\sum_{j=1}^\infty
j\sup_{i\in\mathcal{S}\setminus\{j\},\,x\in\mathbb{R}^m}q_{ij}(x)
<\infty.$$

\noindent $\mathbf{B_4})$\ \ There exists $V\in
C^{1,2}([0,\infty)\times\mathbb{R}^m;\mathbb{R}_+)$ satisfying the
following conditions:

  \noindent (i)
$$ \lim_{|x|\rightarrow\infty}
\left[\inf_{t\in[0,\infty)}V(t,x)\right]=\infty,
$$
and
$$
\lim_{n\rightarrow\infty}\sup_{|x|>n,\,i\in\mathcal{S},\,t\in[0,\infty)}{\cal
L_i}V(t,x) =-\infty.
$$

 \noindent (ii) For any $\rho\ge 1$, there exist $q_{\rho}\in
 B_{b,loc}(\mathbb{R}_+)$ and $W_{\rho}(t,x)\in B_{b,loc}([0,\infty)\times\mathbb{R}^m;\mathbb{R}^m)$ satisfying for each $n\in\mathbb{N}$
there exists $R_n\in L_{loc}^1([0,\infty);\mathbb{R}_+)$ such that
for any $t\in [0,\infty)$ and $x,y\in\mathbb{R}^m$ with
$|x|\vee|y|\leq n$,
$$
   |W_{{\rho}}(t,x)-W_{{\rho}}(t,y)|^2\leq
  R_n(t)|x-y|^2,
  $$
and
$$
\sup_{x\in\mathbb{R}^m,\,i\in\mathcal{S},\,t\in[0,\infty)}{\cal
L_i}V^{\bullet {\rho}}(t,x)<\infty,
$$
\begin{eqnarray*}
V^{\bullet {\rho}}(t,x)\le\langle W_{{\rho}},\nabla_x V^{\bullet
{\rho}}(t,x)\rangle,\ \ \forall t\ge 0,x\in\mathbb{R}^m.
\end{eqnarray*}
Obviously, $\mathbf{Q_3}$ implies $\mathbf{Q_1}$ and
$\mathbf{B_4}$ implies $\mathbf{B_1}$-$\mathbf{B_3}$.

\begin{thm}\label{lv}
Suppose that assumptions $\mathbf{A_2}$, $\mathbf{A_3}$,
$\mathbf{B_4}$,  $\mathbf{Q_2}$ and $\mathbf{Q_3}$ hold. Then,

(i) The hybrid system given by (\ref{aa}) and (\ref{dts}) has a unique $\theta$-periodic solution $(X(t),\Lambda(t))$.

(ii) The Markovian transition semigroup $\{P_{s,t}\}$ of $(X(t),\Lambda(t))$ is strongly Feller and irreducible.

(iii) Let $\mu_s(A)=\mathbb{P}((X(t),\Lambda(t))\in A)$ for $A\in
\mathcal{B}(\mathbb{R}^m\times\mathcal{S})$ and $s\ge 0$. Then,
for any
$s\ge 0$ and $\varphi\in L^2(\mathbb{R}^m\times\mathcal{S};\mu_s)$, we have
$$
\lim_{n\rightarrow\infty}
\frac{1}{n}\sum_{i=1}^{n}P_{s,s+i\theta}\varphi
=\int_{\mathbb{R}^m\times\mathcal{S}}\varphi
d\mu_s\ \
{\rm in}\ \ L^2(\mathbb{R}^m\times\mathcal{S};\mu_s).
$$
\end{thm}
{\bf Proof.} Let $(X(t),\Lambda(t))$ be the unique strong solution to the hybrid system given by (\ref{aa}) and (\ref{dts}) with initial value $(x,i)\in \mathbb{R}^m\times\mathcal{S}$. Define
$$\widetilde{V}(t,x,i)=V(t,x)+i,\ \ t\ge 0,x\in\mathbb{R}^m,i\in{\cal S},$$
where $V$ is given in the assumption $\mathbf{B_4}$. Note that
$$\mathcal{A}\widetilde{V}(t,x,i)=\mathcal{L}_iV(t,x)+\sum_{j\in\mathcal{S}\setminus\{i\}}q_{ij}(x)(j-i).$$
Then, we obtain by assumptions $\mathbf{B_4}$ and $\mathbf{Q_3}$
that
\begin{eqnarray}\label{bq1}
 \lim_{|x|+i\rightarrow\infty}
\left[\inf_{t\in[0,\infty)}\widetilde{V}(t,x,i)\right]=\infty,
\end{eqnarray}
\begin{equation}\label{b3}
\lim_{n\rightarrow\infty}\sup_{|x|+i>n,\,t\in[0,\infty)}{\cal
A}\widetilde{V}(t,x,i) =-\infty,
\end{equation}
and
\begin{equation}\label{bco}
\sup_{x\in\mathbb{R}^m,\,i\in\mathcal{S},\,t\in[0,\infty)}{\cal
A}\widetilde{V}(t,x,i)<\infty.
\end{equation}

For $n\in \mathbb{N}$, we define the stopping time $T_n$ by
$$T_n=\inf\{t\in[0,\infty):|X(t)|\vee \Lambda(t)\geq n\}.$$
For $t\ge 0$, by It\^{o}'s formula, we get
\begin{equation}\label{H43}
\mathbb{E}[\widetilde{V}(t\wedge T_n,X(t\wedge T_n),
\Lambda(t\wedge T_n))]
=\mathbb{E}[\widetilde{V}(0,X(0),\Lambda(0))]+\mathbb{E}\left[\int_{0}^{t\wedge
T_n} {\cal A}\widetilde{V}(u,X(u),\Lambda(u))\mathrm{d}u\right].
\end{equation}

Define
$$
A_n:=-\sup_{|y|+k>n,\,t\in[0,\infty)}{\cal A}\widetilde{V}(t,y,k).
$$
By (\ref{b3}), we get
\begin{equation}\label{H33}
\lim_{n\rightarrow\infty}A_n=+\infty.
\end{equation}
We have
$$
{\cal A}\widetilde{V}(u,X(u),\Lambda(u))\leq -1_{\{|X(u)|+k\geq
n\}}A_n+\sup_{|y|+k<n,\,u\in[0,\infty)}{\cal
A}\widetilde{V}(u,y,k).
$$
Then,  there exist positive constants $c_1$ and $c_2$ such that for large $n$,
\begin{equation}\label{H23}
\mathbb{E}\left[\int_{0}^{t\wedge T_n}
1_{\{|X(u)|+k\geq n\}}\mathrm{d}u\right]\leq \frac{c_{1}t+c_2}{ A_n}.
\end{equation}
Denote $B_n=\{(y,k)\in\mathbb{R}^m\times\mathcal{S}:|y|+k< n\}$ and $B_n^c=\{(y,k)\in\mathbb{R}^m\times\mathcal{S}:|y|+k\ge n\}$.
Letting $n\rightarrow\infty$ in (\ref{H23}),
we obtain by (\ref{H33}) that
\begin{eqnarray}\label{eq1}
\lim_{n\rightarrow\infty}\varlimsup_{T\rightarrow\infty}
\frac{1}{T}\int_{0}^{T}P(0,(x,i),u,B_n^c)\mathrm{d}u=0.
\end{eqnarray}

By (\ref{bco}), there exists
$\lambda>0$ such that
\begin{equation}\label{H44}{\cal A}\widetilde{V}(t,x,i)\le \lambda,\ \ \forall\,t\geq0,\,(x,i)\in \mathbb{R}^m\times\mathcal{S}.
\end{equation}
By (\ref{H43}) and  (\ref{H44}), we get
$$\mathbb{E}[\widetilde{V}(t,X(t),\Lambda(t))]\leq \lambda t+\widetilde{V}(0,x,i).$$
Together with Chebyshev's inequality, this implies that
\begin{equation}\label{H55}\mathbb{P}(0,(x,i),t,B_n^c)\leq \frac{\lambda t+\widetilde{V}(0,x,i)}{\inf_{|y|+k>n,\,t\in[0,\infty)}\widetilde{V}(t,y,k)}.\end{equation}
By (\ref{bq1}) and (\ref{H55}), we find that there exists a
sequence of positive integers $\gamma_n\uparrow\infty$ such that
\begin{eqnarray}\label{eq2}
\lim_{n\rightarrow\infty}\left\{\sup_{(x,i)\in
B_{\gamma_n},\,t\in(0,\theta)}\mathbb{P}(0,(x,i),t,B_n^c)\right\}=0.
\end{eqnarray}

Similar to \cite[Lemma 3.10]{fazxc}, we can show that the
transition semigroup $\{P_{s,t}\}$ of the Markov process
$(X(t),\Lambda(t))$ is $\theta$-periodic. Combining this with
(\ref{eq1}), (\ref{eq2}) and \cite[Theorem 3.2 and Remark
3.1]{a8}, we conclude that the hybrid  system given by (\ref{aa})
and (\ref{dts}) has a $\theta$-periodic solution. Here we would
like to call the reader's attention to a missing condition in
\cite[Theorem 3.2]{a8}, which was pointed out by Hu and Xu
recently. According to Hu and Xu \cite[Theorem 2.1 and Remark
A.1]{a12}, \cite[Theorem 3.2 and Remark 3.1]{a8} holds under the
additional assumption that $\{P_{s,t}\}$ is a Feller semigroup. By
Theorem \ref{thm-4}, $\{P_{s,t}\}$ is a strong Feller semigroup
and hence we can apply \cite[Theorem 3.2 and Remark 3.1]{a8} to
show that the hybrid system given by (\ref{aa}) and (\ref{dts})
has a $\theta$-periodic solution.

By using the same argument, we can show that \cite[Lemmas 3.12 and
3.13]{fazxc} hold with the state  space $\mathbb{R}^m$ replaced by
$\mathbb{R}^m\times\mathcal{S}$. Hence the uniqueness of the
$\theta$-periodic solution is a direct consequence of Theorem
\ref{thm-4} and Theorem \ref{thm-6}. Finally, the last assertion
of the theorem can be proved by following the same argument of the
proof of \cite[Lemma 3.13]{fazxc}. \hfill\fbox

\section{Examples}\label{sec4}\setcounter{equation}{0}

In this section, we use two examples to illustrate our main results.
\begin{exa} (Stochastic Lorenz equation  with regime switching)

We consider the following stochastic Lorenz equation (\cite{lo})
with regime switching:
\begin{eqnarray}\label{Lorenz}
dX_1(t) &=& [-\alpha(t,\Lambda(t))X_1(t)+\alpha(t,\Lambda(t))X_2(t)]dt+\sum_{j=1}^{3}\sigma_{1j}(t,X(t),\Lambda(t))dB_{j}(t)\nonumber\\
&& +\int_{\{|u|<1\}}H_1(t,X(t-),\Lambda(t-),u)\widetilde{N}(dt,du)+\int_{\{|u|\ge 1\}}G_1(t,X(t-),\Lambda(t-),u)N(dt,du),\nonumber\\
  dX_2(t) &=& [\mu(t,\Lambda(t))X_1(t)-X_2(t)-X_1(t)X_3(t)]dt+\sum_{j=1}^3\sigma_{2j}(t,X(t),\Lambda(t))dB_{j}(t)\nonumber\\
 &&+\int_{\{|u|<1\}}H_2(t,X(t-),\Lambda(t-),u)\widetilde{N}(dt,du)
 +\int_{\{|u|\ge 1\}}G_2(t,X(t-),\Lambda(t-),u)N(dt,du),\nonumber\\
 dX_3(t) &=& [-\beta(t,\Lambda(t))X_3(t)+X_1(t)X_2(t)]dt+\sum_{j=1}^3\sigma_{3j}(t,X(t),\Lambda(t))dB_{j}(t)\nonumber\\
 &&+\int_{\{|u|<1\}}H_3(t,X(t-),\Lambda(t-),u)\widetilde{N}(dt,du)+\int_{\{|u|\ge 1\}}G_3(t,X(t-),\Lambda(t-),u)N(dt,du),\nonumber\\
\end{eqnarray}
where $\Lambda(t)$ takes values in $\mathcal{S}=\{1,2,\dots\}$ and
is generated by $Q=(q_{ij}(x))$ that satisfies assumptions
$\mathbf{Q_2}$ and $\mathbf{Q_3}$,
$\alpha(t,i),\beta(t,i),\mu(t,i):[0,\infty)
\times\mathcal{S}\rightarrow\mathbb{R}_{+}$,
$\sigma(t,x,i):[0,\infty)\times\mathbb{R}^3\times\mathcal{S}\rightarrow
\mathbb{R}^{3\times 3}$,
$H(t,x,i,u):[0,\infty)\times\mathbb{R}^3\times\mathcal{S}\times\mathbb{R}^l
\rightarrow\mathbb{R}^3$ and
$G(t,x,i,u):[0,\infty)\times\mathbb{R}^3\times\mathcal{S}\times\mathbb{R}^l
\rightarrow\mathbb{R}^3$ are all Borel measurable and periodic
w.r.t. $t$ with period $\theta$.

We assume that $\sigma(t,x,i)$ and $H(t,x,i,u)$ satisfy
(\ref{qww}) and (\ref{Q}), and there exist $\gamma>0$ and a
continuously differentiable periodic function $a(t)$ with period
$\theta$ such that for any $t\in[0,\theta)$, $i\in\mathcal{S}$,
\begin{eqnarray}\label{KKL}
\gamma<\alpha(t,i),\ \beta(t,i), \
\mu(t,i)\leq a(t).
\end{eqnarray}
Moreover, for any $\varepsilon>0$ there exists $c_{\varepsilon}>0$ such that
for any $i\in\mathcal{S}$,
\begin{eqnarray*}
|\sigma(t,x,i)|^2+\int_{\{|u|<1\}}|H(t,x,i,u)|^2\nu(du)+\int_{\{|u|\ge
1\}}|G(t,x,i,u)|^2\nu(du) \leq \varepsilon|x|^{2}+c_{\varepsilon}.
\end{eqnarray*}

Define
$$V(t,x)=x_1^2+x_2^2+(x_3-2a(t))^2.$$
We will show that $V$ satisfies the assumption $\mathbf{B_4}$.
Without loss of generality, we assume that ${\rho}=1$. The
verification for the case that $\rho>1$ is completely similar. Let
$$W(t,x)=\frac{1}{2}(x_1,x_2,x_3-2a(t)).$$
Then,
\begin{eqnarray*}
   V(t,x)=\langle W,V_x(t,x)\rangle.
\end{eqnarray*}

For $\varepsilon>0$ and $i\in\mathcal{S}$, we obtain by (\ref{KKL}) that
\begin{eqnarray*}
{\cal L_i}V(t,x) &\leq& 4(2a(t)-x_3)a'(t)-2\alpha(t,i)x_1^2-2x_2^2-2\beta(t,i) [x_3^2-2a(t)x_3]+|\sigma(t,x,i)|^2\\
&&+\sum_{j=1}^3\int_{\{|u|<1\}}|H_j(t,x,i,u)|^2\nu(du)
  +\sum_{j=1}^2\int_{\{|u|\ge 1\}}[(x_j+G_j(t,x,i,u))^2-x_j^2]\nu(du)\\
&&+\int_{\{|u|\ge 1\}}[(x_3+G_3(t,x,u)-2a(t))^2-(x_3-2a(t))^2]\nu(du)\\
  &\leq&4(2a(t)+|x_3|)|a'(t)|{-2\left[\gamma x_1^2+x_2^2+\gamma x_3^2-2(a(t))^2|x_3|\right]+\varepsilon|x|^{2}+c_{\varepsilon}}\\
 &&+2[|x_1|+|x_2|+|x_3|+2a(t)]\left[\nu(\{|u|\ge1\})\right]^{1/2}(\varepsilon|x|^{2}+c_{\varepsilon})^{1/2}.
\end{eqnarray*}
Then,
$$\lim_{n\rightarrow\infty}\sup_{|x|>n,\,i\in\mathcal{S},\,t\in[0,\infty)}{\cal
L_i}V(t,x) =-\infty, \ \
\sup_{x\in\mathbb{R}^m,\,i\in\mathcal{S},\,t\in[0,\infty)}{\cal
L_i}V(t,x)<\infty.
$$
Thus, the assumption $\mathbf{B_4}$ is satisfied. Therefore, the
SDE (\ref{Lorenz}) has a unique $\theta$-periodic solution
$(X(t),\Lambda(t))$ and assertions (ii) and (iii) of Theorem
\ref{lv} hold.
\end{exa}

\begin{exa} (Stochastic equation of the lemniscate of Bernoulli with regime switching)

In this example, we consider the stochastic equation of the
lemniscate of Bernoulli, which generalizes \cite[Example 3.20]{CD}
to the non-autonomous case with L\'evy noise and regime-switching.

For $x=(x_1,x_2)\in \mathbb{R}^2$, define
$$I(x)=(x_1^2+x_2^2)^2-4(x_1^2-x_2^2).
$$ Let
$$
{\cal V}(I)=\frac{I^2}{2(1+I^2)^{3/4}},\ \ \ \ {\cal H}(I)=\frac{I}{(1+I^2)^{3/8}}.
$$
Consider the vector field
$$
b(x)=-\left[{\cal V}_x(I)+\left(\frac{\partial {\cal H}(I)}{\partial x_2},-\frac{\partial {\cal H}(I)}{\partial x_1}\right)^T\right].
$$
We have
\begin{eqnarray*}
&&\frac{d{\cal V}(I)}{dI}=\frac{I(I^2+4)}{4(1+I^2)^{7/4}},\ \ \ \ \frac{d{\cal H}(I)}{dI}=\frac{I^2+4}{4(1+I^2)^{11/4}},\label{lemirr}\\
&&\frac{\partial I}{\partial x_1}=4x_1(x_1^2+x_2^2)-8x_1,\ \ \ \ \frac{\partial I}{\partial x_2}=4x_2(x_1^2+x_2^2)+8x_2.\label{lemi1}
\end{eqnarray*}
Define
\begin{eqnarray*}
f(I)=\frac{d{\cal V}(I)}{dI},\ \ \ \ g(I)=\frac{d{\cal H}(I)}{dI}.
\end{eqnarray*}
Then,
\begin{eqnarray*}{\cal V}_x(I)=\frac{d{\cal V}(I)}{dI}\left(\frac{\partial I}{\partial x_1},\frac{\partial I}{\partial x_2}\right)^T,
\end{eqnarray*}
and
\begin{eqnarray*}
&&b_1(x)=-f(I)(4x_1(x_1^2+x_2^2)-8x_1)-g(I)(4x_2(x_1^2+x_2^2)+8x_2),\nonumber\\
&&b_2(x)=-f(I)(4x_2(x_1^2+x_2^2)+8x_2)-g(I)(-4x_1(x_1^2+x_2^2)+8x_1).
\end{eqnarray*}

We consider the following SDEs:
\begin{eqnarray}\label{ssll}
dX_1(t) &=& b_1(X(t))dt+\sigma_{11}(t,X(t),\Lambda(t))dB_1(t)+\sigma_{12}(t,X(t),\Lambda(t))dB_2(t)\nonumber\\
 &&+\int_{\{|u|<1\}}H_1(t,X(t-),\Lambda(t-),u)\widetilde{N}(dt,du)+\int_{\{|u|\ge 1\}}G_1(t,X(t-),\Lambda(t-),u)N(dt,du),\nonumber\\
  dX_2(t) &=& b_2(X(t))dt+\sigma_{21}(t,X(t),\Lambda(t))dB_1(t)+\sigma_{22}(t,X(t),\Lambda(t))dB_2(t)\nonumber\\
 &&+\int_{\{|u|<1\}}H_2(t,X(t-),\Lambda(t-),u)\widetilde{N}(dt,du)
 +\int_{\{|u|\ge 1\}}G_2(t,X(t-),\Lambda(t-),u)N(dt,du),\nonumber\\
\end{eqnarray}
where $\Lambda(t)$ takes values in $\mathcal{S}=\{1,2.\ldots\}$ and is generated by $Q(x)=(q_{ij}(x))$ with
$$q_{ij}(x)=\frac{1\wedge|x|}{j^{2+\delta}} \ \ {\rm for\ some}\ \delta>0,$$
$\sigma(t,x,i):[0,\infty)\times\mathbb{R}^2\times\mathcal{S}\rightarrow
\mathbb{R}^{2\times 2}$,
$H(t,x,i,u):[0,\infty)\times\mathbb{R}^2\times\mathcal{S}\times\mathbb{R}^l
\rightarrow\mathbb{R}^2$ and
$G(t,x,i,u):[0,\infty)\times\mathbb{R}^2\times\mathcal{S}\times\mathbb{R}^l
\rightarrow\mathbb{R}^2$ are all Borel measurable and periodic
w.r.t. $t$ with period $\theta$. We  assume that $\sigma(t,x,i)$
and $H(t,x,i,u)$ satisfy (\ref{qww}) and (\ref{Q}). Moreover, for
any $\varepsilon>0$, there exists $c_{\varepsilon}>0$ such that
for any $i\in\mathcal{S}$, {\small \begin{eqnarray*}
&&|\sigma(t,x,i)|^2+\int_{\{|u|<1\}}|H(t,x,i,u)|^2\nu(du)+|x|\int_{\{|u|\ge 1\}}|G(t,x,i,u)|\nu(du)\\
&&\ \ \ \ \ \ \ \ \ \ \ \ \ \ \ \ \ \ \ \ \ \ +\int_{\{|u|\ge 1\}}|G(t,x,i,u)|^2\nu(du)
\leq \varepsilon|x|^{2}+c_{\varepsilon}.
\end{eqnarray*}}

Define
$$V(t,x)={\cal V}(I(x)).$$
Following the argument of \cite[Example 4.3]{fazxc}, we can check
that the assumption $\mathbf{B_4}$ is satisfied. Therefore, the
SDE (\ref{ssll}) has a unique $\theta$-periodic solution
$(X(t),\Lambda(t))$ and assertions (ii) and (iii) of Theorem
\ref{lv} hold.
\end{exa}

\bigskip



\begin{thebibliography}{1234}






\bibitem{CD} Chen L., Dong Z., Jiang J. and Zhai J.,
On limiting behavior of stationary measures for stochastic
evolution systems with small noise intensity. Sci. China Math. 63
(2020) 1463-1504.

\bibitem{ACFZ0} Chen X., Chen Z., Tran K. and Yin G.,
Properties of switching jump diffusions: maximum principles and Harnack inequalities. Bernoulli 25 (2019) 1045-1075.

\bibitem{ACFZ} Chen X., Chen Z., Tran K. and Yin G.,
Recurrence and ergodicity for a class of regime-switching
jump diffusions. Appl. Math. Optim. 80 (2019) 415-445.







\bibitem{fazxc} Guo X. and Sun W.,
Periodic solutions of stochastic differential
equations driven by L\'{e}vy noises.
https://arxiv.org/pdf/1905.04843v2.pdf.

\bibitem{a12} Hu H. and Xu L., Existence and uniqueness theorems for periodic
Markov process and applications to stochastic functional differential
equations. J. Math. Anal. Appl. 466 (2018) 896-926.




\bibitem{a8} Khasminskii R.Z., Stochastic Stability of
Differential Equations. Springer-Verlag, Second Edition, 2012.

\bibitem{lo} Lorenz E.N., Deterministic nonperiodic flow.
J. Atmos. Sci. 20 (1963) 130-141.


\bibitem{fam} Mao X. and Yuan C., Stochastic Differential
Equations with Markovian Switching. Imperial College Press, London, 2006.

\bibitem{AAQ} Meyn S.P. and Tweedie R.L.,
Markov Chains and Stochastic Stability.
Communications and Control Engineering Series.
Springer-Verlag London, 1993.

\bibitem{NY1} Nguyen D.H. and Yin G., Modeling and analysis of switching diffusion systems: past-dependent switching with a countable state space. SIAM J. Control Optim. 54 (2016) 2450-2477.

\bibitem{NY2} Nguyen D.H. and Yin G., Recurrence and ergodicity of switching diffusions with past-dependent switching having a countable state space. Potent. Anal. 48 (2018) 405-435.


\bibitem{AAE} Shao J., Strong solutions and strong Feller
properties for regime-switching diffusion processes in an
infinite state space. SIAM J. Control Optim. 4 (2015) 2462-2479.

\bibitem{ACFA} Xi F., Asymptotic properties of jump-diffusion processes with
state-dependent switching.  Stoch. Proc. Appl. 119 (2009) 2198-2221.

\bibitem{BBC} Xi F. and Yin G., Jump-diffusions with state-dependent
switching: existence and uniqueness, Feller property, linearization,
and uniform ergodicity. Sci. China Math. 12 (2011) 2651-2667.

\bibitem{zvaqw} Xi F., Yin G. and Zhu C.,
Regime-switching jump diffusions with
non-Lipschitz coefficients and countably many
switching states: existence and uniqueness,
Feller, and strong Feller properties. Modeling, Stochastic Control, Optimization, and Applications. 571-599, 2019.


\bibitem{FFRTA} Xi F. and Zhu C.,
On Feller and strong Feller properties and expeniential ergodicity
of regime-switching jump diffusion processes with countable regimes.
SIAM J. Control. Optim. 55 (2017) 1789-1818.









\bibitem{YZ} Yin G. and Zhu C., Hybrid Switching Diffusions: Properties and Applications. Stochastic Medelling and Applied Probability. Vol. 63. Springer, New York, 2010.



\bibitem{k1} Zhang X., Wang K. and Li D.,
Stochastic periodic solutions of stochastic differential
equations driven by  L\'{e}vy process.
J. Math. Anal. Appl. 430 (2015) 231-242.








\end{thebibliography}
\end{document}